\documentclass[a4paper,11pt]{article}
\usepackage[plainpages=false]{hyperref}
\usepackage{amsfonts,latexsym,rawfonts,amsmath,amssymb,amsthm,mathrsfs}
\usepackage{amsmath,amssymb,amsfonts,latexsym,lscape,rawfonts}

\usepackage[all]{xy}
\usepackage{eufrak}
\usepackage{makeidx}         

\usepackage{graphicx,psfrag}
\usepackage{array,tabularx}

\usepackage{setspace}

\newtheorem{thm}{Theorem}[section]

\newtheorem{lem}[thm]{Lemma}

\theoremstyle{remark}
\newtheorem{rmk}[thm]{Remark}

\theoremstyle{definition}

\def \R {\mathbb R}
\def \D {\mathbb D}

\def \D {\mathcal D}
\def \Q {\mathcal Q}

\def \U {\mathcal U}
\def \V {\mathcal V}

\def \W {\mathcal W}

\def \NB {\mathcal N}

\def \L {\mathcal L}

\def \g {\mathfrak g}
\def \p {\partial}
\def \gd {\dot{g}(t)}
\def \Riem {\mathcal R}
\def \Wb {\widetilde{\W}}

\DeclareMathOperator{\Hess}{Hess}

\DeclareMathOperator{\dVol}{dVol}

\DeclareMathOperator{\Vol}{Vol}

\DeclareMathOperator{\Ker}{Ker}

\DeclareMathOperator{\Image}{Im}

\DeclareMathOperator{\Diff}{Diff}

\begin{document}
\title{On the K\"ahler-Ricci flow near a K\"ahler-Einstein metric}
\author{Song Sun and Yuanqi Wang}
\date{}
\maketitle{}

\begin{abstract}
On a Fano manifold, we prove that the K\"ahler-Ricci flow starting from a K\"ahler metric in the anti-canonical class which is sufficiently close to a K\"ahler-Einstein metric
must converge in a polynomial rate to a K\"ahler-Einstein metric. The convergence can not happen in general if we study the flow on the level of K\"ahler potentials. Instead we exploit the interpretation of the Ricci flow as  the gradient flow of Perelman's $\mu$ functional. This involves modifying the Ricci flow by a canonical family of gauges. In particular, the complex structure of the limit could be different in general.  The main technical ingredient is a Lojasiewicz type inequality for Perelman's $\mu$ functional near a critical point.
\end{abstract}

\section{Introduction}
In this article we shall prove the following
\begin{thm}\label{main theorem}
Let $(M, \omega, J, g)$ be a K\"ahler-Einstein manifold of
Einstein constant 1, i.e. $Ric(g)=g$. Then there is a $C^{k+10,
\gamma}$($k\gg1$) tensor neighborhood $\NB$ of $(\omega, J, g)$,
such that the following holds. For any K\"ahler structure
$(\omega', J', g')$ which lies in $\NB$ and satisfies $\omega'\in
2\pi c_1(M;J')$, we denote by $(\omega(t), J(t), g(t))$ the
(normalized) K\"ahler-Ricci flow starting from $g'$
\begin{equation}\label{Kahler-Ricci flow}
\left\{
                                                               \begin{array}{lll}
\frac{\p g(t)}{\p t}=-Ric(g(t))+g(t),   & t\geq 0\\
J(t)= J' &t\geq 0\\
          g(0)=g'.    & \\

                                                               \end{array}
                                                             \right.
\end{equation}
Then the flow exists for all $t\in [0,\infty)$, and  there is an isotopy of diffeomorphisms $f_t$ such that
$f_t^*(\omega(t), J(t), g(t))$ converges in $C^{k, \gamma}$ to a
K\"ahler-Einstein metric $(\omega_{\infty}, J_{\infty},
g_{\infty})$ in a polynomial rate, i.e. there are constants $C>0$,
$\alpha>0$, such that
\begin{equation}
||f_t^*g(t)-g_{\infty}||_{C^{k, \gamma}}+||f_t^*J(t)-J_{\infty}||_{C^{k, \gamma}}\leq C\cdot
(t+1)^{-\alpha},
\end{equation}
where the $C^{k, \gamma}$   topology is taken with respect to the
initial K\"ahler-Einstein metric $(\omega, J,g)$. Moreover, if $J$
is adjacent to $J'$ in the sense that there is a sequence of
diffeomorphisms $\phi_i$ such that $\phi_i^*J'\rightarrow J$(see
\cite{CS}), then the limit $(\omega_{\infty}, J_{\infty},
g_{\infty})$ is isomorphic to $(\omega, J, g)$ under the action of
the diffeomorphism group.
\end{thm}

\begin{rmk}
We should mention that two other approaches on the stability of
Ricci flow on Fano manifolds have been previously announced by
Arezzo-La Nave(see \cite{AL}) and by Tian-Zhu(see \cite{TZ3}). Our
method here is different and more direct, and Theorem 1.1 follows
from a stability lemma for general modified Ricci flow
(\ref{stability lemma}).
\end{rmk}

R. Hamilton(\cite{Ha}) introduced the \emph{Ricci flow} as a way
to produce Einstein metrics. He proved the short time existence
using his own generalization of the Nash-Moser inverse function
theorem. Later De Turck reproved the short time existence by  a
gauge fixing trick. The long time existence for general Ricci flow
fails. H-D. Cao(\cite{Cao}) first studied the Ricci flow on a
K\"ahler manifold, called the \emph{K\"ahler-Ricci flow}. The
K\"ahler condition is preserved under the flow, and the equation
can be reduced to a scalar equation on the K\"ahler potential.
Using Yau's estimates for complex Monge-Amp\`ere equations,  Cao
proved the long time existence of the flow in any K\"ahler class
proportional to the canonical class, and proved convergence in the
case of non-positive first Chern class. In the case of positive
first Chern class, the K\"ahler-Ricci flow in the anti-canonical
class $2\pi c_1(M)$ takes the normalized form as in
(\ref{Kahler-Ricci flow}) with the complex structure fixed. In
this case the flow does not necessarily converge to a
K\"ahler-Einstein metric, due to the obstructions to the existence
of K\"ahler-Einstein metrics on Fano manifolds. In an unpublished
paper, Perelman proved that the scalar curvature and diameter are
always uniformly bounded along the K\"ahler-Ricci flow, and he
announced that if there exists a K\"ahler-Einstein metric, then
the flow starting from any K\"ahler metric in the anti-canonical
class of the same complex structure converges to a
K\"ahler-Einstein metric(Note this was previously proved by
Chen-Tian(\cite{CT}) when the initial metric is assumed to have
positive bisectional curvature). In \cite{TZ1}, using Perelman's
estimates, an alternative proof of this convergence was given. The
case we consider is a stability type theorem: the K\"ahler-Ricci
flow initiating from a K\"ahler metric near a K\"ahler-Einstein
metric(possibly with a different complex structure) always
converges(modulo diffeomorphisms) polynomially fast to a
K\"ahler-Einstein metric. Using the gradient interpretation of
Ricci flow by Perelman(\cite{Pe}), Tian-Zhu(\cite{TZ2}) previously
proved such a stability theorem when the K\"ahler-Einstein metric
is assumed to be isolated. Our theorem is based on a more close
investigation of this gradient nature, and is motivated by the
case of Calabi flow studied in \cite{CS}. We believe similar idea
may be helpful in the study of other geometric flows. Note the
gauging diffeomorphisms must diverge, if there is no
K\"ahler-Einstein metric compatible with the complex structure
$J'$. There are also various studies on the stability of the Ricci
flow near a Ricci flat metric, for instance Sesum(\cite{Se1})
proved exponential convergence of the Ricci flow near a Ricci flat
metric which is linearly stable and integrable. For more details
on this topic, readers are referred to \cite{Se1} and the
references therein. In our case, the convergence is exponential
precisely when the complex structure $J'$ itself admits a
K\"ahler-Einstein metric, for the exponential convergence of $g(t)$ would imply the exponential decay of the Ricci potential of $g(t)$, which guarantees the exponential convergence of the K\"ahler potential. In general, there are obstructions(c.f. \cite{Do},  \cite{Sz}) to
deform K\"ahler-Einstein metrics, so the
exponential convergence fails.  \\

The idea of the proof is as follows. Perelman discovered that the
Ricci flow is essentially a gradient flow.  Let $\Riem(M)$ be the
space of all Riemannian metrics on $M$. There is an appropriate
Riemannian metric defined on $\Riem(M)$, which is invariant under
the action of the diffeomorphism group $\Diff(M)$. The Ricci flow
direction differs from the gradient of the so-called $\mu$
functional only by an infinitesimal action of $\Diff(M)$. Both the
Ricci flow and the $\mu$ functional are $\Diff(M)$-invariant, and
thus live on $\Riem(M)/\Diff(M)$. Then on this quotient, the Ricci
flow in the normalized form as in (\ref{Kahler-Ricci flow}) is
exactly the gradient flow of the $\mu$ functional. If $[g_0]$ is a
local maximum of $\mu$ on $\Riem(M)/\Diff(M)$, then we ask whether
the Ricci flow initiating from a nearby point would always stay
close to $[g_0]$ and converge to a maximum point of $\mu$ at
infinity. If we were in finite dimension, then by Lojasiewicz's
fundamental structure theorem for real analytic functions(\cite{Lo}), the
flow would converge polynomially fast to a unique limit at
infinity provided the functional is real-analytic. The limit is
also a local maximum, but possibly be different from the one we
start with if there is a non-trivial moduli. Note in the case when
the maximum is non-degenerate in the sense of Morse-Bott, then
indeed we have exponential convergence. But if the maximum is
degenerate, then we can not expect exponential convergence, and
the rate of convergence depends on how bad the degeneracy is.  For
more details about the finite dimensional setting, the readers are
referred to \cite{CS}.  Our problem has an infinite dimensional
nature, but as in \cite{Si}, the real difficulty is still finite
dimensional. Heuristically, one can decompose the tangent space at
$[g_0]$ into the direct sum of two parts. One part is infinite
dimensional but on which the Hessian of $\mu$ is non-degenerate,
and the other part is finite dimensional on which we can apply
Lojasiewicz's inequality. We need to check the functional $\mu$ is
real-analytic near a K\"ahler-Einstein metric. The last problem is
that our K\"ahler-Einstein metric may not be a local maximum of
$\mu$ among all variations, but is so among all K\"ahler metrics
in the same real cohomology class. This is enough for our purpose,
thanks to the fact that such a
subset is preserved under the Ricci flow. \\

This paper is organized as follows. In section 2, we set up some
notations and definitions. In section 3, we prove a Lojasiewicz
type inequality for the $\mu$ functional(Lemma \ref{Lojasiewicz
for mu}), and prove a general stability theorem for the modified
Ricci flow(Lemma \ref{stability lemma}). In section 4, we
prove Theorem \ref{main theorem}, using results in section 3.\\

\textbf{{Acknowledgements:}} Both authors would like to thank
Professor Xiuxiong Chen for constant support. We also thank Prof.
C. Arezzo and Prof. X-H. Zhu for their interest in this paper and
for sending us their preprints(see \cite{AL}, \cite{TZ3}). The
first author would also like to thank the department of
Mathematics in Stony Brook for its hospitality during the year
2009-2010. The first author is partially supported by an Research
Assistantship in an NSF grant. We thank the anonymous referee for a careful reading and pointing out several typos in the first version.

\section{Preliminaries}
Suppose $M$ is an $n$-dimensional smooth compact manifold. Denote
by $\Riem(M)$ the space of all Riemannian metrics on $M$. This is
an open convex subset of the linear space of all smooth sections
of symmetric 2-tensors on $M$. Later we will assign either the
$C^{k,\gamma}$ topology or $L^2_k$ topology on $\Riem(M)$. Right
now we simply take the $C^{\infty}$ topology. For any smooth
measure $dm$ on $M$, we denote by
$C^{\infty}_0(M;dm;\R)$($C^{k,\gamma}_0(M;dm;\R)$) the space of
$C^{\infty}$($C^{k,\gamma}$) real-valued functions $f$ on $M$
satisfying
$$\int_M e^{-f}dm=1.$$ Here and afterwards, when we use a norm without mentioning the metric, it is always meant to be the norm taken with respect to the fixed metric $g_0$. Given a Riemannian metric $g$, we define
Perelman's functional
\begin{equation}\label{def of W}
\W_g(f)=\int_M [\frac{1}{2}(|\nabla f|^2+R(g))+f]e^{-f}\dVol_g,
\end{equation}
for any $f\in C^{\infty}_0(M;\dVol_g;\R)$.
 We
will denote this functional either by $\W_g(f)$ or $\W(g, f)$,
where in the latter case we emphasize $\W$ as a function depending
on both the metric $g$ and the function $f$. Then by a
straightforward calculation the first variation of $\W$ is given
by(c.f. \cite{Pe})
\begin{lem}\label{1st variation for W}
\begin{eqnarray}\label{1st variation of W}
\delta\W(h, v)=&&\int_M[-\frac{1}{2}\langle Ric(g)+\Hess f-(\Delta
f-\frac{1}{2}|\nabla f|^2+f+\frac{1}{2}R(g))g,
h\rangle_{g}\nonumber\\&&-(\Delta f-\frac{1}{2}|\nabla
f|^2+f+\frac{1}{2}R(g)-1)v]e^{-f}\dVol_g.
\end{eqnarray}
\end{lem}

For a given Riemannian metric $g$, we define Perelman's $\mu$
functional to be
\begin{equation}
\mu(g)=\inf_{f\in C^{\infty}_0(M;\dVol_g;\R)} \W_g(f)
\end{equation}
By a minimizing procedure(see \cite{Ro}) the infimum is always
achievable by a function $f$ which possesses the same regularity
as the metric $g$. From equations (\ref{def of W}) and (\ref{1st variation of W}) we see
that $f$ satisfies the non-linear equation
\begin{equation}\label{minimizer equation for W}
\Delta f-\frac{1}{2}|\nabla f|^2+f+\frac{1}{2}R(g)=\mu(g).
\end{equation}

We call a metric $g\in \Riem(M)$ \emph{regular} if there is a
neighborhood $\U$ of $g$, such that for any $g'\in\U$, the
minimizer of $\W_{g'}$ is unique and depends real-analytically on
$g'$. The following lemma will be proved in the next section.

\begin{lem}\label{Einstein is regular}
A normalized shrinking gradient Ricci soliton(i.e. $Ric(g)+\Hess
f=g$) is regular.
\end{lem}

Now we assume $g$ is regular, then it is easy to see from
(\ref{1st variation of W}) the first variation of $\mu$ is given
by
\begin{equation}\label{1st variation of mu}
\delta\mu(h)=-\frac{1}{2}\int_M \langle Ric(g)+\Hess f-g,
h\rangle_ge^{-f}\dVol_g,
\end{equation}
where $f$ is the minimizer of $\W_g$. If we endow $\Riem(M)$ with
the Riemannian metric
\begin{equation}\label{Riemannian metric}
(h_1, h_2)_g:=\frac{1}{2}\int_M \langle h_1, h_2\rangle_ge^{-f}\dVol_g,
\end{equation}
then we see that
\begin{equation}\label{gradient of mu}
\nabla\mu(g)=-(Ric(g)+\Hess f-g)
\end{equation}
So the critical points of $\mu$ are exactly normalized shrinking
gradient Ricci solitons. The gradient flow of $\mu$ functional is
\begin{equation}\label{modified Ricci flow}
\left\{
                                                               \begin{array}{ll}
\frac{\p g(t)}{\p t}=-Ric(g(t))+g(t)-\Hess_tf(t),   & t\geq 0\\
          g(0)=g.    & \\

                                                               \end{array}
                                                             \right.
\end{equation}
We shall call this flow  the \emph{``modified" Ricci flow}.  We
have
\begin{lem} \label{two flows}Up to an isotopy of diffeomorphisms, the gradient flow (\ref{modified Ricci flow}) is
equivalent to the normalized Ricci flow
\begin{equation}\label{normalized Ricci flow}
\left\{
                                                               \begin{array}{ll}
\frac{\p g(t)}{\p t}=-Ric(g(t))+g(t),   & t\geq 0\\
          g(0)=g,    &\\

                                                               \end{array}
                                                             \right.
\end{equation}
if we assume $g(t)$ is regular for all time as long as the flow
exists.
\end{lem}
\begin{proof} This is because under our hypothesis the function $f(t)$ depends smoothly on $t$, and then we
can translate between the two flows by the isotopy of
diffeomorphisms generated by $\nabla_t f(t)$.
\end{proof}

\section{A stability lemma for the Ricci flow}
 We
shall prove in the end of this section the following Lojasiewicz
type inequality.

\begin{lem}  \label{Lojasiewicz for mu}Let $g_0$ be a normalized shrinking gradient Ricci soliton. Then there
is a $C^{k,\gamma}$($k\gg1$) neighborhood $\U$ of $g_0$ in
$\Riem(M)$, and constants $C>0$, and $\alpha\in [\frac{1}{2}, 1)$,
such that for any $g\in \U$, we have
\begin{equation}\label{Lojasiewicz inequality}
||\nabla \mu(g)||_g\geq C\cdot |\mu(g_0)-\mu(g)|^{\alpha},
\end{equation}
where $||\cdot||_g$ denotes the $L^2$ norm taken with respect to $g$.
\end{lem}
For convenience, we denote
$$\V^{k, \gamma}_{\delta}=\{g\in \Riem(M)|||g-g_0||_{C^{k,\gamma}_{g_0}}\leq \delta\}.$$

\begin{lem}\label{stability lemma}
Suppose $g_0$ is a normalized shrinking gradient Ricci soliton.
Then there exist $\delta_2>\delta_1>0$, such that for any $g\in
\V^{k+10, \gamma}_{\delta_1}$, the modified Ricci flow
\begin{equation}
\left\{
                                                               \begin{array}{ll}
\frac{\p g(t)}{\p t}=-Ric(g(t))+g(t)-\Hess_tf(t),   & t\geq 0\\
          g(0)=g ,    & t=0\\

                                                               \end{array}
                                                             \right.
\end{equation}
starting from any $g$ satisfies  $g(t)\in\V^{k,
\gamma}_{\delta_2}$ as long as  $\mu(g(t))\leq\mu(g_0)$. In
particular, if we know a priori that $\mu(g(t))\leq \mu(g_0)$ for
all $t$, then the flow $g(t)$ exists globally for all time $t$,
and converges in $C^{k, \gamma}_{g_0}$ to a limit $g_{\infty}$ which is
also a shrinking gradient Ricci soliton with
$\mu(g_{\infty})=\mu(g_0)$. The convergence is in a polynomial
rate.
\end{lem}

\begin{proof} The proof is by Lojasiewicz arguments(see \cite{CS}). For convenience, we include the details here.
Choose $k$ large, $\gamma\in(0,1)$ and $\delta>0$ small such that all metrics in
$\V^{k,\gamma}_{\delta}$ are regular with (\ref{Lojasiewicz
inequality}) holds for uniform constants $C$ and $\alpha$. By the short time stability(see Lemma \ref{Modified Ricci flow's finite time stability} below) we know that there are $T_0>0$ and $\delta_1>0$ such that for all $g\in \V^{k+10, \gamma}_ {\delta_1}$ the modified Ricci flow $g(t)$ starting from $g$ exists and lies in $\V^{k, \gamma}_{\delta/4}$ for $t\in[0,T_0]$. So it suffices to prove that there exists
$0<\delta_1<\delta_2<\delta$ such that for any modified Ricci flow
solution $g(t)$ with $g(0)\in\V^{k+10, \gamma}_{\delta_1}$, if
$g(t)\in \V^{k,\gamma}_{\delta}$ and $\mu(g(t))\leq\mu(g_0)=0$ for
$t\in[0, T)$($T>T_0$), then $g(T)\in \V^{k, \gamma}_{\delta_2}$(see
figure \ref{fig: stability}).
\begin{figure}
 \begin{center}
  \psfrag{A}[c][c]{$\V^{k+10, \gamma}_{\delta_1}$}
  \psfrag{B}[c][c]{$\V^{k, \gamma}_{\delta_2}$}
  \psfrag{C}[c][c]{$\V^{k, \gamma}_{\delta}$}
  \psfrag{D}[c][c]{$g(t)$}
  \includegraphics[width=0.8 \columnwidth]{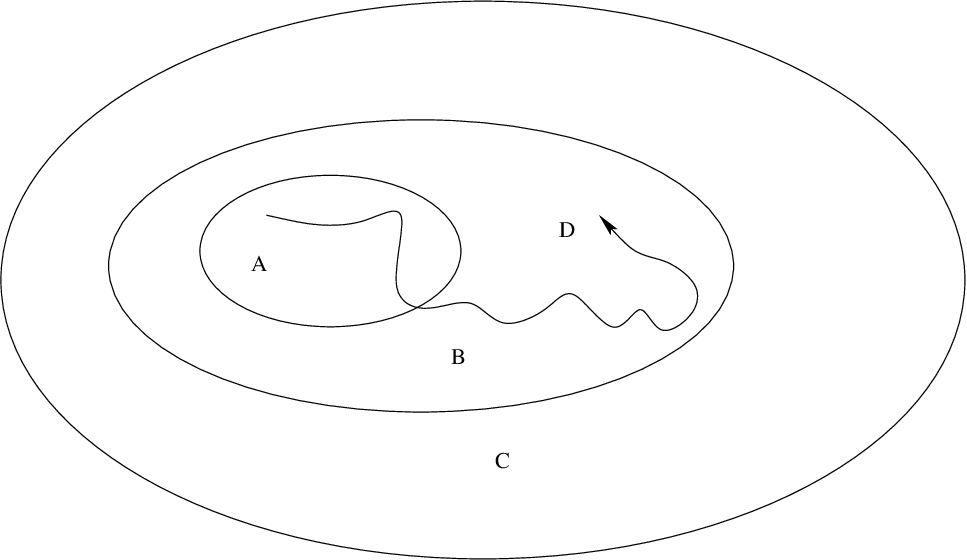}
  \caption{The flow $g(t)$ can never exit $\V^{k, \gamma}_{\delta_2}$}
  \label{fig: stability}
 \end{center}
 \end{figure}
 By the fundamental curvature estimates of
Hamilton,   for any integer $l\geq 1$, and $t\in[T_0,T)$, we have
$$||Rm(g(t))||_{C^{l, \gamma}_t}\leq C(l), $$
where the norm is taken with respect to the metric $g(t)$. Here and below we adopt the notation that $C(\star)$ is a constant only depending on $\star$ and the already fixed data $g_0$, $k$, $\gamma$,  $\delta$ and $T_0$, but the precise value of $C(\star)$ could vary from line to line and is not relevant for our purpose. Since all
metrics in $\V^{k,\gamma}_{\delta}$ are regular, this gives
\begin{equation}
||f(t)||_{C^{l, \gamma}_t}\leq C(l).
\end{equation}
Fix $\beta\in (2-\frac{1}{\alpha}, 1)$, where $\alpha$ is the constant in Lemma \ref{Lojasiewicz for mu}.   Then by interpolation
inequalities for tensors,  for any integer $p\geq1$ there is an
$N(p)$(independent of $t\geq 1$), such that
\begin{eqnarray*}
||\gd||_{L^2_p(t)}&\leq&
C(p)\cdot||\gd||^{\beta}_{L^2(t)}\cdot||Ric(g(t))-g(t)+\Hess_t
f(t)||^{1-\beta}_{L^2_{N(p)}(t)} \\&\leq& C(p)\cdot
||\gd||^{\beta}_{L^2(t)}.
\end{eqnarray*}
Here again $||\cdot||_{L^2_p(t)}$ denotes the Sobolev norm taken with respect to the metric $g(t)$. Since
\begin{eqnarray}
&&\frac{d}{dt}[\mu(g_0)-\mu(g(t))]^{1-(2-\beta)\alpha}\nonumber\\&=&
-[1-(2-\beta)\alpha]\cdot[\mu(g_0)-\mu(g(t))]^{-(2-\beta)\alpha}\cdot||\nabla\mu(g(t))||^2\nonumber\\&\leq&
-C(\beta)\cdot || \nabla\mu(g(t))||^{\beta}_{L^2(t)},
\end{eqnarray} by integration we get
\begin{eqnarray}\label{length estimate}
&&\int_{T_0}^T||\gd||^{\beta}_{L^2(t)}dt\nonumber \\&\leq&
C(\beta)^{-1}\cdot [(\mu(g_0)-\mu(g(T_0)))^{1-(2-\beta)\alpha}-(\mu(g_0)-\mu(g(T)))^{1-(2-\beta)\alpha}]\nonumber\\&\leq&
 C(\beta)^{-1}\cdot [\mu(g_0)-\mu(g(0))]^{1-(2-\beta)\alpha}.
\end{eqnarray}
Therefore
$$\int_{T_0}^T||\gd||_{L^2_p(t)}dt\leq C(p)\cdot \epsilon(\delta_1),$$
where $\epsilon(\star)$ denotes a quantity that goes to zero as $\star$ goes to zero and the data $g_0$,  $k$, $\gamma$ and $\delta$  are fixed.
Since the Sobolev constant is uniformly bounded for metrics in $\V^{k,
\gamma}_{\delta}$, we obtain for any $l\geq1$,
\begin{equation}\int_{T_0}^T||\gd||_{C^{l,\gamma}_t}dt\leq C(l)\cdot
\epsilon(\delta_1),\end{equation} and
$$||g(T)-g(T_0)||_{C^{k,\gamma}_t}\leq
\int_{T_0}^T||\gd||_{C^{k,\gamma}_t}dt=\epsilon(\delta_1).$$ Since the
$C^{k,\gamma}$ norm defined by metrics in $\V^{k,
\gamma}_{\delta}$ are equivalent to each other, with bounds only depending on $\delta$, we obtain
$$||g(T)-g(T_0)||_{C^{k,\gamma}_{g_0}}=\epsilon(\delta_1).$$
Since $$||g(T_0)-g_0||_{C^{k, \gamma}_{g_0}}\leq \delta/4,  $$ we obtain 
$$||g(T)-g_0||_{C^{k, \gamma}_{g_0}}=\epsilon(\delta_1)+\delta/4.$$
Now choose $\delta_2=\frac{\delta}{2}$, and choose $\delta_1$ so that
$\epsilon(\delta_1)\leq \delta/4$, then the first part of the
lemma is proved.

Now we assume $\mu(g(t))\leq \mu(g_0)$ for all $t$, then $g(t)$
exists for all time. Indeed, $g(t)$ can never exit
$\V^{k, \gamma}_{\delta_2}$.
Since \begin{eqnarray*}
&&\frac{d}{dt}[\mu(g_0)-\mu(g(t))]^{1-2\alpha}\nonumber\\&=&
-(1-2\alpha)\cdot[\mu(g_0)-\mu(g(t))]^{-2\alpha}\cdot||\nabla\mu(g(t))||^2\nonumber
\\&\geq&C(2\alpha-1)\cdot[\mu(g_0)-\mu(g(t))]^{-2\alpha}\cdot[\mu(g_0)-\mu(g(t))]^{2\alpha}
\\&\geq&
C(2\alpha-1).
\end{eqnarray*}
This implies
\begin{equation*}
\mu(g_0)-\mu(g(t))\leq (C(2\alpha-1))^{-\frac{1}{2\alpha-1}}\cdot t^{-\frac{1}{2\alpha-1}}.
\end{equation*}
Since  $|\mu(g)-\mu(g_0)|=\epsilon(\delta_1)$ hods for $g\in \V^{k+10,\gamma}_{\delta_1}$, we have $\mu(g_0)-\mu(g(t))\le \mu(g_0)-\mu(g(0))\leq \epsilon(\delta_1)$ for all $t$.   Together with (\ref{decay estimate}) we obtain
\begin{equation} \label{decay estimate}
\mu(g_0)-\mu(g(t))\leq C(\beta) \cdot (t+1)^{-\frac{1}{2\alpha-1}}.
\end{equation}
Then for any $t_2\geq t_1$, the same argument as before shows
\begin{eqnarray*}
||g(t_1)-g(t_2)||_{C^{k, \gamma}_{g_0}}&\leq&
\int_{t_1}^{t_2}||\dot{g}(t)||_{C^{k, \gamma}_{g_0}}dt\\&\leq&
C(\beta)\cdot [\mu(g_0)-\mu(g_{t_1})]^{1-(2-\beta)\alpha}\\& \leq&
C(\beta)\cdot (t_1+1)^{-\frac{1-(2-\beta)\alpha}{2\alpha-1}}.
\end{eqnarray*}
Hence the flow $g(t)$ converges uniformly in $C^{k, \gamma}_{g_0}$ to a
limit $g_{\infty}$, and let $t_2\rightarrow\infty$ we obtain
polynomial decay rate. By (\ref{length estimate}) and (\ref{decay
estimate}), we see that $\nabla\mu(g_{\infty})=0$, and
$\mu(g_{\infty})=\mu(g_0)$. So $g_{\infty}$ is a shrinking
gradient Ricci soliton.
\end{proof}

\begin{lem}\label{Modified Ricci flow's finite time stability}
Let $(M, g_0)$ be a normalized shrinking Ricci soliton with $Ric(g_0)+\Hess f_0=g_0$. Then for any $k\gg1$ and $\epsilon>0$, there exists a $T_0>0$, and $\delta>0$  such that for any $g$ with
$||g-g_0||_{C^{k+10,\gamma}_{g_0}}\leq \delta$,  the modified Ricci flow starting  from $g$
exists for $t\in [0,T_0]$. Moreover, $||g(t)-g_0||_{C^{k, \gamma}_{g_0}}\leq \epsilon$ for all $t\in [0,T_0]$.
\end{lem}
This is a standard result that follows from the De Turck trick for short time existence of Ricci flow and the standard theory on quasilinear parabolic equations. The possible loss of derivatives is due to the fact that we have to  transform from De Turck Ricci flow to Ricci flow, and then to the modified Ricci flow by diffeomorphisms.  This is a technical point and is not relevant for our main purpose. \\

Now we start to prove Lemma \ref{Einstein is regular}. We first show the
minimizer of $\W_g$ is unique if $g$ is a shrinking gradient Ricci
soliton.

\begin{lem}\label{minimizer unique for soliton} If $g$ is a shrinking gradient Ricci soliton: $$Ric(g)+\Hess
f=g$$ with $\int_M e^{-f}\dVol_g=1$,  then the minimizer of $\W_g$
is unique and is equal to $f$. In particular, if $Ric(g)=g$, then
the minimizer of $\W_g$ is $\log \Vol_g(M)$.
\end{lem}
\begin{proof} Let $\phi_t$ be the one-parameter group of diffeomorphisms generated by $\nabla f$, and denote $g(t)=
\phi_t^*g$. Then $$\frac{\p}{\p t}g(t)=\Hess_t
\phi_t^*f=g(t)-Ric(g(t)).$$ Suppose $v$ is any minimizer of
$\W_g$,
 Now solve the backward heat equation
\begin{equation}
\left\{
                                                               \begin{array}{ll}
\frac{\p v(t)}{\p t}=\frac{1}{2}[n-R(g(t))-\Delta_t v(t)]+|\nabla_t v(t)|^2,   & t\in [0, \tau]\\
          v(\tau)=\phi_{\tau}^*v    & \\

                                                               \end{array}
                                                             \right.
\end{equation}
in a small time interval $[0, \tau]$. Let $\psi_t$ be the isotopy
of diffeomorphisms generated by the (time-dependent) vector fields
$-\nabla_t v(t)$. Denote  $\tilde{g}(t)=\psi_t^*g(t)$, and
$\tilde{v}(t)=\psi_t^*v(t)$, then
$$\frac{\p \tilde{g}(t)}{\p t}=\tilde{g}(t)-Ric(\tilde{g}(t))-\Hess_t \tilde{v}(t), $$
and
$$\frac{d\tilde{v}}{dt}=\frac{1}{2}[n-R(\tilde{g}(t))-\Delta_t\tilde{v}(t)].$$
Then by (\ref{1st variation of W}) we have
$$\frac{\p}{\p t}\W(\tilde{g}(t), \tilde{v}(t))=\frac{1}{2}\int_M|\tilde{g}(t)-Ric(\tilde{g}(t))-\Hess_t\tilde{v}(t))|^2
e^{-\tilde{v}(t)}\dVol_t\geq 0.$$ However, since all
$\tilde{g}(t)$ are in the same $\Diff(M)$ orbit, we have
$\mu(\tilde{g}(t))\equiv \mu(g)$. Since $\phi_{\tau}^*g$ is a
minimizer of $\W_{g(\tau)}$, we have
$$\W(\tilde{g}(t), \tilde{v}(t))\equiv \W(g(\tau), \phi_{\tau}^*v)=\W(g, v). $$
Hence $$g-Ric(g)-\Hess v=0, $$ and so $v=f$.
 \end{proof}

 In general, for any function
$f$, we denote by $d^{*f}$ and $\nabla^{*f}$ the adjoint of the
$d$ and $\nabla$  with respect to the measure $e^{-f}\dVol_g$.
Perelman(\cite{Pe}) observed the following Bochner formula
relating the twisted Hodge Laplacian
$\Delta_H^{f}=d^{*f}d+dd^{*f}$ and rough Laplacian
$\Delta^f=-\nabla^{*f}\nabla$:
\begin{equation}
-\Delta^f\xi=\Delta_H^{f}\xi-Ric^f(\xi),
\end{equation}
where $\xi$ is any one-form and $Ric^f=Ric+\Hess f$. Standard Weitzenb\"ock formula gives

\begin{lem}\label{eigenvalue estimate}
Suppose $Ric(g_0)+\Hess f_0=g_0$, then the first nonzero
eigenvalue of $-\Delta^{f_0}$ acting on functions is strictly bigger than one.
\end{lem}

 \textbf{Proof of lemma \ref{Einstein is regular}}.
 Suppose $Ric(g_0)+\Hess f_0=g_0$. Define
$$L: C^{k,\gamma}(\Riem(M))\times C^{k, \gamma}(M;\R)\rightarrow C^{k-4,\gamma}(M;\R),$$
which sends $(g, f)$ to
$$\Delta_g^f(\Delta_gf+f-\frac 1 2|\nabla_gf|^2+\frac 1 2 R(g))+\int_Me^{-f}\dVol_g-1.$$
Then $L$ is a real-analytic map between Banach manifolds. We have
$L(g_0, f_0)=0$, and the differential of $L$ at $(g_0, f_0)$ has
its second component equal to
\begin{equation}
DL_{(g_0,f_0)}(0, f)=\Delta_{g_0}^{f_0}(\Delta_{g_0}^{f_0}f+f)
-\int_Mf\dVol_{g_0}.
\end{equation}
This is an isomorphism by Lemma \ref{eigenvalue estimate}. Thus by
the real-analytic version of the  implicit function theorem for
Banach manifolds (c.f \cite{Ko}), there is
 a $C^{k, \gamma}$
neighborhood $\V_1$ of $g_0$ in $\Riem(M)$  and a real-analytic
map $P: \V_1\rightarrow C^{k,\gamma}(M;\R)$ such that for any
$g\in \V_1$, $P(g)\in C^{k,\gamma}_0(M;\dVol_g;\R)$, and
$$L(g,P(g))=0.$$ Moreover, there is a $\delta>0$ such that if
$L(g, f)=0$ for some $g\in \V_1$ and  $||f-f_0||_{C^{k,\gamma}}\leq
\delta$, then $f=P(g)$. Now $P(g)\in C^{k,
\gamma}_0(M;\dVol_g;\R)$ satisfies the equation
 $$\Delta_gP(g)+P(g)-\frac 1 2 |\nabla_gP(g)|_g^2+\frac 1 2 R(g)=\text{constant}.$$ In particular,
$P(g)$ is a critical point of $\W_g$. Now we claim that we can
choose $\V_2\subset\V_1$, such that for any $g\in \V_2$, the
minimizer for $\W_g$ is unique and is equal to $P(g)$. Suppose
this were false, then there would be a sequence $g_i\in \V_1$, and
a minimizer $f_i$ of $\W_{g_i}$, such that $g_i\rightarrow g_0$ in
$C^{k,\gamma}$ topology, and $f_i\neq P(g_i)$ for any $i$. Let
$u_i=e^{-\frac{f_i}{2}}$, then $u_i$ is a minimizer of the
functional
$$\Wb_{g_i}(u)=\frac{1}{2}\int_M 4|\nabla_i u|_i^2+R(g_i)u^2-2u^2\log u^2 \dVol_{g_i}$$
under the constraint
$$\int_Mu_i^2\dVol_{g_i}=1.$$ Here and after  an operator with subscript $i$ denotes the operator corresponding to the metric $g_i$, and $|\cdot|_i$ denotes the pointwise norm taken with respect to $g_i$.
 Moreover, we have $\Wb_{g_i}(u_i)=\mu(g_i)$ and
\begin{equation}\label{elliptic equation}
-4\Delta_i u_i+R(g_i)u_i-4u_i\log u_i-2\mu(g_i)u_i=0.
\end{equation}
By definition,  Perelman's $\mu$ functional is upper
semi-continuous, i.e.
$$\overline{\lim_{i\rightarrow\infty}}\mu(g_i)\leq \mu(g_0).$$
Thus for all $i$ we have
$$\int_M 4|\nabla_i u_i|_i^2+R(g_i)u_i^2-4u_i^2\log u_i\dVol_{g_i}=2\mu(g_i)\leq
C.$$ Here $C$ is a constant which does not depend on $i$, but
might vary from line to line.  We claim $||u_i||_{L^{\infty}}\leq C$. This is standard from \cite{Ro}. For the convenience of the reads, we include the proof here.  We assume $n>2$, as the case $n=2$ can be done similarly.
 By Jensen's inequality,
$$\int_M2u_i^2\log u_i\dVol_{g_i}=\frac{n-2}{2}\int_Mu_i^2\log u_i^{\frac{4}{n-2}}\dVol_{g_i}\leq \frac{n-2}{2}\log \int_M
u_i^{\frac{2n}{n-2}}\dVol_{g_i}.$$
For any positive real number $a$ we have a positive number $b(a)$ such that $$\log x\leq ax^{\frac{n-2}{2n}}+b(a).$$
Then we have $$\int_Mu_i^2\log u^2_i d Vol_{g_i}\leq \frac{(n-2)a}{2}(\int_M  u^{\frac{2n}{n-2}}_i d Vol_{g_i})^{\frac{n-2}{2n}}+\frac{n-2}{2}b(a).$$
Then using $\mu(g_i)\leq C$ we have
$$\int_M|\nabla_i u_i|^2_i d Vol_{g_i}\leq \frac{(n-2)a}{4}(\int_M  u^{\frac{2n}{n-2}}_i d Vol_{g_i})^{\frac{n-2}{2n}}+\frac{n-2}{4}b(a)+C.$$
 Since $g_i\in\V_1$, we have a
uniform bound of the Sobolev constant, so
$$(\int_M
u_i^{\frac{2n}{n-2}}\dVol_{g_i})^{\frac{n-2}{2n}}\leq  C\cdot (||\nabla_i
u_i||_{L^2_{g_i}}^2+||u_i||_{L^2_{g_i}}^2), $$
and then
$$\int_M|\nabla_i u_i|_i^2 d Vol_{g_i}\leq \frac{(n-2)a}{4}C\int_M|\nabla_i u_i|_i^2 d Vol_{g_i}+\frac{n-2}{4}[b(a)+Ca]+C.$$
Choose $a$ to be sufficiently small to make $\frac{(n-2)a}{4}C\leq \frac{1}{2}$ we have that:
$$\int_M|\nabla_i u_i|_i^2 d Vol_{g_i}\leq C.$$
So $$||u_i||_{L^{\frac{2n}{n-2}}_{g_i}}\leq C, $$
and $$||u_i\log u_i||_{L^{\frac{2n}{n-2+\epsilon}}_{g_i}}\leq C, $$
for any $\epsilon>0$ small. Since $g_i\in \V_1$, we have uniform elliptic estimate for $\Delta_i$. Thus by equation (\ref{elliptic equation}), and Soboleve embedding we see
$$||u_i||_{L^{\frac{2n}{n-6+\epsilon}}_{g_i}}\leq C.$$
Thus again $$||u_i\log u_i||_{L^{\frac{2n}{n-6+\epsilon_1}}_{g_i}}\leq C, $$
for $\epsilon_1>\epsilon$ sufficiently small.
Keep on iterating and after at most $n/4+1$ steps we arrive at  the claim
$$||u_i||_{L^{\infty}}\leq C. $$ Then for any $\beta\in(0,1)$,
$$||u_i||_{C^{1,\beta}_{g_i}}\leq C.$$
Note by assumption the $C^{k, \gamma}$ norm defined by all the
metrics $g_i$ are equivalent to that defined by the metric $g_0$.
By passing to a subsequence we can assume $u_i$ converges to a
limit $u_{\infty}$ in $C^{1, \beta}_{g_0}$. Since $u_i$ is a positive
minimizer of $\Wb_{g_i}$, $u_{\infty}$ is a non-negative minimizer
of $\Wb_{g_0}$. By the strong maximum principle in \cite{Ro},
$u_{\infty}$ is strictly positive. So we obtain a uniform positive
lower bound for $u_i$. Then we get a uniform $C^{1,\beta}_{g_0}$ bound
on $f_i$. Now by applying elliptic regularity for $f_i$, we get
$$||f_i||_{C^{k,\gamma}_{g_0}}\leq C.$$ So by passing to a subsequence,
$f_i$ converges in weak $C^{k, \gamma}_{g_0}$ topology to $f_{\infty}$.
Then it is easy to see that $f_{\infty}=-2\log u_{\infty}$ is a
minimizer for $\W_{g_0}$. By Lemma \ref{minimizer unique for
soliton} we obtain $f_{\infty}=f_0$.  Note again implicit function
theorem ensures that for $i $ sufficiently large any critical
point of $\W_{g_i}$ which is $C^{k, \gamma'}_{g_0}$($\gamma'$ slightly
smaller than $\gamma$) close to $f_0$ must be $P(g_i)$. Thus
$f_i=P(g_i)$, and we arrive at a contradiction.
$\square$\\

\textbf{Proof of Lemma \ref{Lojasiewicz for mu}} Denote by
$\Riem'(M)$ the space of regular Riemannian metrics on M. $\Riem'(M)$ is endowed with the Riemannian metric given by
equation (\ref{Riemannian metric}). We put the $L^2_k$ topology on
$\Riem'(M)$(for convenience we do not use the H\"older norm here)
and denote the completion by $\Riem'_k(M)$, so that $\Riem'_k(M)$
becomes a Banach manifold. Near a point $g_0$, $\Riem'_k(M)$ can be
identified with an open set in the Banach space $\Gamma_k(g_0)$,
which is the space of $L^2_k$ sections of the bundle of symmetric
two-tensors $h$. Again as before, when we use a norm without mentioning the metric, we mean the norm taken with respect to $g_0$. This gives rise to a coordinate chart for
$\Riem'_k(M)$. We shall then identify any $h\in \Riem'_k(M)$ close
to $g_0$ with its ``coordinate" $h\in\Gamma_k(g_0)$, by abuse of
notation.

Let $\Diff_{k+1}(M)$ be the group of $L^2_{k+1}$ diffeomorphisms
of $M$.   It acts on $\Riem'(M)$, preserving the Riemannian metric.
Notice that both $\mu$ and $||\nabla\mu||$ are invariant under the
action of $\Diff_{k+1}(M)$. Here $||\nabla\mu||(g)=||\nabla\mu(g)||_g$ is the $L^2$ norm of $\nabla\mu(g)$ with respect to $g$, which appears on the left hand side of  equation \ref{Lojasiewicz inequality}.
In a neighborhood of the identity
map, $\Diff_{k+1}(M)$ is modelled on the linear space $\g_{k+1}$
of $L^2_{k+1}$ vector fields $X$ on $M$. At $g_0$, the tangent
space to the $\Diff_{k+1}(M)$ orbit of $g_0$ is given by the image
of
$$\L_k: \g_{k+1} \rightarrow \Gamma_k(g_0); X\mapsto \nabla^s X, $$
where $\nabla^s X$ is the symmetrization of $\nabla X$. The normal
space to $\Image \L_k$ with respect to the Riemannian metric is
given by $\Ker \L_k^*=\ker \nabla^{*f_0}$, which consists of
$L^2_k$ divergence free symmetric 2-tensors. By the standard slice
theorem(see for example \cite{Eb}), there are neighborhoods
$\U\subset \V$ of $g_0$ in $\Riem_k(M)$ for any $g\in \U$, there
is a $\phi\in \Diff_{k+1}(M)$ close to identity such that
$\phi^*g=h\in \V\cap \Ker \L_k^*$. Here we implicitly made use of
the identification mentioned above. Thus it suffices to prove
inequality (\ref{Lojasiewicz inequality}) for $g$ in $\Q=\V\cap
\Ker \L_k^*\subset \Gamma_k(g_0)$. We still denote by $\mu$ its
restriction to $\Q$, and $\nabla_{\Q}\mu$ the gradient of $\mu$ on
$\Q$, with respect to the $L^2$ metric(defined by $g_0$). Then $\nabla_{\Q}\mu\in \Ker \L_{k-2}^*$ and we
have for any $g\in \Q$
$$||\nabla\mu(g)||_g\geq C\cdot ||\nabla_{\Q}\mu(g)||_{L^2}.$$
So it suffices to prove
\begin{lem}
 There exists a neighborhood $\NB\subset \Q$ of $g_0$ and constant
$C>0$, $\alpha>0$,  such that for any $g\in\NB$, we have
\begin{equation}\label{Lojasiewicz for restriction}
||\nabla_{\Q}\mu(g)||_{L^2}\geq C\cdot |\mu(g)-\mu(g_0)|^{\alpha}.
\end{equation}
\end{lem}
\begin{proof}
We follow the same pattern of arguments as in \cite{CS}. By
(\ref{gradient of mu}), $0$ is a critical point of $\mu$. By a
direct computation, the Hessian of $\mu$ at $0$ (viewed as an
operator from $\Ker \L_k^*$ to $\Ker \L_{k-2}^*$) is given by:
\begin{eqnarray}
H_0(h)&=&\frac{1}{2}\Delta^{f_0} h+Rm(g_0)\circ h,
\end{eqnarray}
which is a ``twisted " Lichnerowicz Laplacian.
  Then $H_0$ is an elliptic differential operator of order 2 on $\Ker \L_k^*$. So it has a
finite dimensional kernel $W_0$ which  consists of smooth
elements, and we have the following decomposition:
$$\Ker \L_k^*=W_0\oplus W_k',$$
where $H_0$ restricts to an invertible operator from $W_k'$ to
$W_{k-2}'$. So there exists a $C>0$, such that for any $\eta'\in
W'$, we have
$$||H_0(\eta')||_{L^2_{k-2}}\geq
C\cdot||\eta'||_{L^2_k}.$$ By the implicit function theorem, there
are small constants $\epsilon_1, \epsilon_2>0$, such that for any
$\eta_0\in W_0$ with $||\eta_0||_{L^2}\leq \epsilon_1$(so
$||\eta_0||_{L^2_{k-2}}$ is also small since $W_0$ is finite
dimensional), there exists a unique element $\eta'=G(\eta_0)\in
W'$ with $||\eta'||_{L^2_k}\leq \epsilon_2$, such that
$\nabla_{\Q}\mu(\eta_0+\eta')\in W_0$. Moreover by construction
the map $G: B_{\epsilon_1}W_0\rightarrow B_{\epsilon_2}W'$ is real
analytic. Now consider the function
$$F: W_0\rightarrow \R; \eta_0\mapsto \mu(\eta_0+G(\eta_0)).$$
This is a real analytic function with $\eta_0=0$ as a critical
point. For any $\eta_0\in W_0$, it is easy to
see that $\nabla F(\eta_0)=\nabla_{\Q}\mu(\eta_0+G(\eta_0))\in W_0$.\\

Now we shall estimate the two sides of the inequality
(\ref{Lojasiewicz for restriction}) separately. For any $\eta\in
W$ with $||\eta||_{L^2_k}\leq \epsilon$, we can write
$\eta=\eta_0+G(\eta_0)+\eta'$, where $\eta_0\in W_0$, $\eta'\in
W'$, and
$$||\eta_0||_{L^2_k}\leq c\cdot ||\eta||_{L^2_k},$$
$$||G(\eta_0)||_{L^2_k}\leq c\cdot ||\eta||_{L^2_k},$$
$$||\eta'||_{L^2_k}\leq c\cdot ||\eta||_{L^2_k}.$$
For the left hand side of $(\ref{Lojasiewicz for restriction})$,
we have:
\begin{eqnarray*}
\nabla_{\Q}\mu(\eta)&=&\nabla_{\Q}\mu(\eta_0+G(\eta_0)+\eta')\\
              &=&\nabla_{\Q}\mu(\eta_0+G(\eta_0))+\int_0^1
              \delta_{\eta'}\nabla_{\Q}\mu(\eta_0+G(\eta_0)+s\eta')ds\\
              &=&\nabla F(\eta_0)+\delta_{\eta'}\nabla_{\Q}\mu(0)+\int_0^1
              (\delta_{\eta'}\nabla_{\Q}\mu(\eta_0+G(\eta_0)+s\eta')-\delta_{\eta'}\nabla_{\Q}
              \mu(0))ds\\
\end{eqnarray*}
The first two terms are $L^2$ orthogonal to each other. For the
second term we have
$$||\delta_{\eta'}\nabla_{\Q}\mu(0)||_{L^2}^2=||H_0(\eta')||_{L^2}^2\geq C\cdot||\eta'||^2_{L^2_2}.$$ For
the last term, we have
$$||\delta_{\eta'}\nabla_{\Q}\mu(\eta_0+G(\eta_0)+s\eta')-\delta_{\eta'}\nabla_{\Q}\mu
              (0)||_{L^2}\leq C\cdot||\eta||_{L^2_k}||\eta'||_{L^2_2}\leq C\cdot\epsilon\cdot ||\eta'||_{L^2_2}. $$
Therefore, we have
\begin{equation}||\nabla_{\Q}\mu(\eta)||_{L^2}^2\geq
|\nabla
F(\eta_0)|_{L^2}^2+C\cdot||\eta'||_{L^2_2}^2.\end{equation}

For the right hand side of $(\ref{Lojasiewicz for restriction})$,
we have
\begin{eqnarray*}
\mu(\eta)&=& \mu(\eta_0+G(\eta_0)+\eta')\\
       &=& \mu(\eta_0+G(\eta_0))+\int_0^1 \nabla_{\Q}\mu(\eta_0+G(\eta_0)+s\eta')\eta'ds\\
       &=& F(\eta_0)+\nabla F(\eta_0)\eta'+\int_0^1\int_0^1
       \delta_{\eta'}\nabla_{\Q}\mu(\eta_0+G(\eta_0)+st\eta')\eta'dtds\\
       &=& F(\eta_0)+H_0(\eta')\eta'+\int_0^1\int_0^1
       (\delta_{\eta'}\nabla_{\Q}\mu(\eta_0+G(\eta_0)+st\eta')-\delta_{\eta'}\nabla_{\Q}\mu(0))\eta'dtds
\end{eqnarray*}
So
\begin{equation}
|\mu(\eta)-\mu(0)|\leq |F(\eta_0)-F(0)|_{L^2}+
C\cdot||\eta'||_{L^2_2}^2.
\end{equation}
Now we apply the  usual Lojasiewicz inequality to $F$, and obtain
that
$$|\nabla F(\eta_0)|_{L^2}\geq C\cdot|F(\eta_0)-F(0)|^{\alpha}, $$
for some $\alpha\in [\frac{1}{2},1)$, and $C>0$. Of course here the constant $C$ may vary from line to line. Together we have proved $(\ref{Lojasiewicz for restriction})$.\\
\end{proof}

\section{Proof of the main theorem}
Now we shall prove Theorem \ref{main theorem}. Suppose $(M,
\omega, J, g)$ is a K\"ahler-Einstein metric of Einstein constant
one. Denote by $\NB^{k+10,\gamma}_{\delta}$ the space of all K\"ahler
structures $(\omega', J', g')$ with $\omega'\in 2\pi\cdot
c_1(M;J')$ and
$$||g-g'||_{C^{k+10,\gamma}_g}+||J-J'||_{C^{k+10,\gamma}_g}\leq\delta.
$$ Then there exists a small $\delta$, such that
any $(\omega', J', g')\in \NB^{k+10,\gamma}_{\delta}$ satisfies
$c_1(M;J')=c_1(M;J)$, and such that
  the inequality
(\ref{Lojasiewicz inequality}) is satisfied for any K\"ahler
metric $(\omega', J', g')\in\NB^{k+10,\gamma}_{\delta}$, by Lemma
\ref{Lojasiewicz for mu}. By definition \begin{equation*}
\mu(g')\leq \W_{g'}(\log \Vol_{g'}(M))=\frac{1}{2}\int_M
R(g')\dVol_{g'}=\frac{n}{2}+\log \Vol_{g'}(M).
\end{equation*}
Moreover, the equality holds if and only if $Ric(g')=g'$(To see this, if $Ric(g')=g'$, then by Lemma \ref{minimizer unique for soliton} we see the equality holds; for the converse, if the equality holds, then the constant function $f=\log \Vol_{g'}(M)$ is a minimizer of $\W_{g'}$, so by equation (\ref{minimizer equation for W}) one sees $R(g')$ is constant. Then since $g'$ is K\"ahler and $\omega'\in 2\pi c_1(M; J')$, we conclude $Ric(g')=g'$.)
Therefore, we have
$$\mu(g')\leq \mu(g).$$
 Since the K\"ahler-Ricci flow preserves the anti-canonical K\"ahler condition,   we
are able to make use of Lemma \ref{stability lemma}, and obtain
constants $0<\delta_1<\delta_2<\delta$ such that  the modified
Ricci flow starting from any $g'\in \NB^{k+10,\gamma}_{\delta_1}$
converges in $C^{k, \gamma}$ to a Ricci soliton $g_{\infty}\in
\NB^{k, \gamma}_{\delta_2}$ polynomially fast with
$\mu(g_{\infty})=\mu(g)$. This latter condition implies that
$g_{\infty}$ is indeed Einstein.
 The complex structure $J(t)$ under the
modified Ricci flow evolves as
$$\frac{d}{dt}J(t)=J(t)\D_{t}f(t), $$
where $\D$ is the Lichnerowicz Laplacian in K\"ahler geometry.
Since
$$\nabla\mu((g(t)))=Ric(g(t))-g(t)+\Hess_t f(t)=[Ric(g(t))-g(t)+\nabla_t\bar{\nabla}_tf(t)]+\D_t f(t)$$ is a
point-wise orthogonal decomposition, we see that
$$||\dot{J}(t)||_{L^2(t)}\leq ||\dot{g}(t)||_{L^2(t)}, $$ thus by
following the argument in the proof of Lemma \ref{stability lemma}, we
see that $J(t)$ also converges in $C^{k, \gamma}$ to a limit
$J_{\infty}$ polynomially fast. Moreover, $J_{\infty}$ and
$g_{\infty}$ are compatible, so the K\"ahler structures
$(\omega(t), J(t), g(t))$ converges in $C^{k, \gamma}$ to the
K\"ahler-Einstein structure $(\omega_{\infty}, J_{\infty},
g_{\infty})$ in $\NB^{k, \gamma}_{\delta_2}$. By Lemma \ref{two
flows}, we conclude the main theorem. The last part follows from
\cite{CS}. $\square$

\begin{rmk} For a K\"ahler-Ricci soliton, since we do not
know whether it is always a maximizer of $\mu$ among nearby
K\"ahler metrics in a fixed real  cohomology class, we can not
conclude stability of K\"ahler-Ricci flow in this case. But
similar arguments using Lemma \ref{stability lemma} can show that
if the K\"ahler-Ricci flow converges by sequence to a
K\"ahler-Ricci soliton in the sense of Cheeger-Gromov, then it
converges uniformly and polynomially fast in the sense of Theorem
\ref{main theorem}. The uniqueness of the limit soliton of a Ricci
flow has been proved by N. Sesum(\cite{Se2}) under the extra
assumption of integrability.
\end{rmk}

Department of Mathematics, Imperial College, London, SW7 2AZ, U.K.

Email: s.sun@imperial.ac.uk.\\

Department of Mathematics, University of California at Santa Barbara, CA 93106, U.S.

Email: wangyuanqi@math.ucsb.edu.

\end{document}